\documentclass[letterpaper,10pt]{amsart}
\usepackage{amsmath,amssymb,amsxtra,amsthm}

\makeatletter
\@namedef{subjclassname@2010}{%
  \textup{2010} Mathematics Subject Classification}
\makeatother

\theoremstyle{plain}
\newtheorem{theorem}{Theorem}[section]

\newtheorem{lemma}[theorem]{Lemma}

\newtheorem{conjecture}[theorem]{Conjecture}

\theoremstyle{definition}

\newcommand{\B}{\mathbb}

\newcommand{\ga}{\alpha}

\newcommand{\eps}{\varepsilon}
\newcommand{\gf}{\varphi}

\newcommand{\gs}{\sigma}

\newcommand{\gz}{\zeta}

\begin{document}
\title{On the multiplicative Erd\H{o}s discrepancy problem}
\author{Michael Coons}
\address{University of Waterloo, Dept.~of Pure Math., Waterloo, ON, N2L 3G1, Canada}
\email{mcoons@math.uwaterloo.ca}
\thanks{The research of M.~Coons is supported by a Fields-Ontario Fellowship and NSERC.}%
\date{\today}

\begin{abstract} As early as the 1930s, P\'al Erd\H{o}s conjectured that: {\em for any multiplicative function $f:\mathbb{N}\to\{-1,1\}$, the partial sums $\sum_{n\leqslant x}f(n)$ are unbounded.} In this paper, after providing a counterexample to this conjecture, we consider completely multiplicative functions $f:\B{N}\to\{-1,1\}$ as well as a class of similar multiplicative functions $f$ satisfying $$\sum_{p\leqslant x}f(p)=c\cdot\frac{x}{\log x}(1+o(1)).$$ We prove that if $c>0$ then the partial sums of $f$ are unbounded, and if $c<0$ then the partial sums of $\mu f$ are unbounded. Extensions of this result are also given.
\end{abstract}

\maketitle

\section{Introduction}

Erd\H{o}s \cite{E1957} asked the following question, sometimes known as the Erd\H{o}s Discrepancy Problem. {\em ``Let $f(n)=\pm 1$ be an arbitrary number theoretic function. Is it true that to every $c$ there is a $d$ and an $m$ for which \begin{equation}\label{E8}\left|\sum_{k=1}^n f(kd)\right|>c\ ?\end{equation} Inequality \eqref{E8} is one of my oldest conjectures.''} (This particular quote is taken from a restatement of the conjecture in \cite[p.78]{E1985}. See also \cite{E1985b} and \cite{EG}.) Erd\H{o}s offered 500 dollars for a proof of this conjecture. Erd\H{o}s \cite[p.293]{E1957} wrote in 1957 that this conjecture is twenty-five years old, placing its origin at least as far back as the early 1930s. In \cite{E1957,E1985,E1985b}, Erd\H{o}s also stated a multiplicative form of his conjecture. 

\begin{conjecture}[Erd\H{o}s]\label{Econj} Let $f(n)=\pm 1$ be a multiplicative function, (i.e., $f(ab)=f(a)f(b)$, when $\gcd(a,b)=1$). Then  \begin{equation}\label{E9}\limsup_{x\to\infty}\left|\sum_{n\leqslant x} f(n)\right|=\infty;\end{equation} that is, the partial sums of $f$ are unbounded.\end{conjecture} 

\noindent Erd\H{o}s added in \cite{E1985} that {\em ``clearly \eqref{E9} would follow from \eqref{E8} but as far as I know \eqref{E9} has never been proved. Incidentally \eqref{E9} was also conjectured by Tchudakoff.''}

Conjecture \ref{Econj} as stated is not true, and while this may be known to others in this field, there seems to be no account of it in the literature. 

For a counterexample, consider the multiplicative function $g$ defined by $g(1)=1$, and on prime powers by \begin{equation}\label{g}g(p^k)=\begin{cases} -1 &\mbox{if $p=2$ and $k\geq 1$}\\ 1 &\mbox{if $p\neq 2$ and $k\geq 1$}.\end{cases}\end{equation} Then $g$ is periodic with period $2$ and for all $n\geq 1$ we have $g(2n)=-1$ and $g(2n-1)=1$. Thus $$\sum_{n\leq x}g(n)=\begin{cases} 1 & \mbox{if $[x]$ is odd}\\ 0 & \mbox{if $[x]$ is even},\end{cases}$$ and so $$\limsup_{x\to\infty}\left|\sum_{n\leqslant x} g(n)\right|=1.$$

It may very well be the case that the function $g$ defined above is the only counterexample to Conjecture \ref{Econj}, but at least at this point, we can say that this is the only known counterexample.

Along with Conjecture \ref{Econj}, Erd\H{o}s \cite{E1957} conjectured a result on the mean values of multiplicative functions. A number--theoretic function $f:\mathbb{N}\to\mathbb{C}$ has a {\em mean value}, denoted $M(f)$, provided the limit \begin{equation}\label{MV} M(f):=\lim_{x\to\infty}\frac{1}{x}\sum_{n\leqslant x}f(n)\end{equation} exists. Erd\H{o}s \cite{E1957,E1985} (among others) conjectured that any multiplicative function taking the values $\pm 1$ has a mean value; this is usually called the Erd\H{o}s--Wintner Conjecture. In 1961, Delange \cite{Del1} characterized those functions with positive mean value, and in 1967, Wirsing \cite{Wir1967} gave a complete solution to this conjecture, as well as the extension to all complex--valued multiplicative functions $f$ satisfying $|f|\leqslant 1$. This was later refined by Hal\'asz \cite{Halasz} in 1968. We state the result here only for those functions with which we are directly concerned.

\begin{theorem}[Delange, Wirsing, Hal\'asz]\label{TMV} Let $f:\mathbb{N}\to\{-1,1\}$ be a multiplicative function. If \begin{equation}\label{1f}\sum_{p\leqslant x}\frac{1-f(p)}{p}\end{equation} is bounded then $M(f)$ exists and is positive, and if \eqref{1f} is unbounded then $M(f)=0$.
\end{theorem}

We note that the ideas of Theorem \ref{TMV} have been generalized by many authors, including Granville and Soundararajan \cite{GS2007, GS2008} and Goldmakher \cite{Leo}. In these works the authors use properties of a generalization of \eqref{1f} to give some new results concerning sums of certain types of Dirichlet characters. The generalization of \eqref{1f} is usually made by considering a special multiplicative function $g$ (e.g., a Dirichlet character) and comparing it to the multiplicative function of interest $f$ (e.g., a Dirichlet character) by means of investigating the asymptotics of $$\sum_{p\leqslant x}\frac{1-\Re(f\overline{g}(p))}{p}.$$ This sum can be thought of as a metric \cite{GS2007}, and in some sense measures how $g$ mimics $f$; this terminology was introduced in \cite{Leo}. 

In contrast to this ``mimicry metric,'' we consider the asymptotics of $$\sum_{p\leqslant x}\frac{c-f(p)}{p}$$ for $c$ not necessarily equal to $1$. By considering sums like like this, we are able to give the following result toward Conjecture \ref{Econj}.

\begin{theorem}\label{Tmain} Let $f:\B{N}\to\{-1,1\}$ be a multiplicative function such that there is some $k\geq 1$ with $f(2^k)=1$. Suppose that for some $c\in[-1,1]$ we have $$\sum_{p\leqslant x}f(p)=c\cdot\frac{x}{\log x}(1+o(1)).$$ If $c>0$ then the partial sums of $f$ are unbounded, and if $c<0$ the partial sums of $\mu f$ are unbounded.
\end{theorem}

Some extensions of this theorem are given in Section \ref{exten}, including some instances of the case $c=0$. In Section \ref{CMF}, we show that this theorem is true for completely multiplicative functions without the assumption that there is some $k\geq 1$ with $f(2^k)=1$.

\section{Completely multiplicative functions}\label{CMF}

If a multiplicative function $f:\B{N}\to\{-1,1\}$ has positive mean value, then clearly the partial sums of $f$ are unbounded; they are asymptotic to $M(f)\cdot x$. The triviality leaves when we consider functions with $M(f)=0$.

\begin{theorem} Let $f:\mathbb{N}\to\{-1,1\}$ be a completely multiplicative function (i.e., $f(ab)=f(a)f(b)$ for all $a,b\in\B{N}$) and suppose that $c\in[-1,1)$. If $$\sum_{p}\frac{c-f(p)}{p}<\infty,$$ then the mean value of $f$ exists and is equal to $0$.
\end{theorem}

\begin{proof} This follows from Theorem \ref{TMV} in a very straightforward way. We need only note that \begin{multline*} \sum_{p\leqslant x}\frac{1-f(p)}{p}=\sum_{n\leqslant x}\frac{1-c+c-f(p)}{p}\\ =(1-c)\sum_{n\leqslant x}\frac{1}{p}+\sum_{n\leqslant x}\frac{c-f(p)}{p}=(1-c)\log\log x+O(1).\qedhere\end{multline*}
\end{proof}

To prove Theorem \ref{Tmain}, we will first prove the result for completely multiplicative functions $f:\B{N}\to\{-1,1\}$. The bulk of the work is taken up by the following lemma.

\begin{lemma}\label{Lcf} Let $f:\mathbb{N}\to\{-1,1\}$ be a completely multiplicative function. Suppose that $c\in[-1,1]$ is nonzero and $$\sum_{p}\frac{c-f(p)}{p}<\infty.$$ If $c>0$ then the partial sums of $f$ are unbounded, and if $c<0$ the partial sums of $\mu f$ are unbounded.
\end{lemma}

\begin{proof} Suppose firstly that $c>0$. To give the desired result, it is enough to show that $$\lim_{x\to\infty} \sum_{n\leqslant x}\frac{f(n)}{n}=\infty.$$ To this end, note that for $\gs>1$ we have \begin{multline}\label{logF} \log F(\gs)=\log \sum_{n\geqslant 1}\frac{f(n)}{n^\gs} =-\sum_p\log\left(1-\frac{f(p)}{p}\right)\\ =\sum_p\sum_{k\geqslant 1}\frac{f(p)^k}{kp^{k\gs}} =\sum_p \frac{f(p)}{p^\gs}+\sum_p\sum_{k\geqslant 2}\frac{f(p)^k}{kp^{k\gs}}=\sum_p \frac{f(p)}{p^\gs}+O(1),\end{multline} where the $O(1)$ term is valid for $\gs>1/2$. Since $$\sum_p\frac{c-f(p)}{p}<\infty,$$ we have that \begin{equation}\label{fposc}\sum_{p\leqslant x}\frac{f(p)}{p}=c\log\log x+O(1).\end{equation} The condition that $c>0$ ensures that $$\lim_{s\to 1^{+}}\sum_{p}\frac{f(p)}{p^\gs}=\infty,$$ and so the divergence of $\log F(\gs)$ at $\gs=1$ occurs because $\lim_{\gs\to 1^+}F(\gs)=\infty$. 

In the light of \eqref{logF} it must be the case that \begin{equation}\label{fninfty}\lim_{x\to\infty}\sum_{n\leqslant x}\frac{f(n)}{n}=\infty.\end{equation} Thus we have that $$\limsup_{x\to\infty} \left|\sum_{n\leqslant x}f(n)\right|=\infty.$$ For if not, there is a real number $M>0$ such that $\left|\sum_{n\leqslant x}f(n)\right|<M,$ and by partial summation, we would then have that \begin{equation*}\sum_{n\leqslant x}\frac{f(n)}{n}=\frac{1}{x}\sum_{n\leqslant x}f(n)+\int_1^x\left(\sum_{n\leqslant t}f(n)\right)\frac{dt}{t^2}=O\left(\int_1^x\frac{dt}{t^2}\right)=O(1),\end{equation*} which contradicts \eqref{fninfty}.

Now suppose that $c<0$. In this case, instead of $F(\gs)$, we consider the function $1/F(\gs)$. Running through the above argument gives \begin{equation}\label{1overF}-\log F(\gs)=-\sum_p \frac{f(p)}{p^\gs}+O(1),\end{equation} where again the $O(1)$ term is valid for $\gs>1/2.$ Similar to the above, using the assumption of the lemma, we have that \begin{equation}\label{fnegc}-\sum_{p\leqslant x}\frac{f(p)}{p}=|c|\log\log x+O(1),\end{equation} which in turn gives, due to \eqref{1overF} that $$\lim_{\gs\to 1^+}\frac{1}{F(\gs)}=\infty.$$ This implies that $$\lim_{x\to\infty}\sum_{n\leqslant x}\frac{\mu(n)f(n)}{n}=\infty,$$ which using a similar argument as the case $c>0$, give that $$\limsup_{x\to\infty}\left|\sum_{n\leqslant x}\mu(n)f(n)\right|=\infty.$$ This completes the proof of the lemma.
\end{proof}

Our proof of the main theorem follows from the similar result for completely multiplicative functions. Using partial summation we have the following theorem.

\begin{theorem}\label{Tcmain} Let $f:\B{N}\to\{-1,1\}$ be a completely multiplicative function. Suppose that for some $c\in[-1,1]$ we have $$\sum_{p\leqslant x}f(p)=c\cdot\frac{x}{\log x}(1+o(1)).$$ If $c>0$ then the partial sums of $f$ are unbounded, and if $c<0$ the partial sums of $\mu f$ are unbounded.
\end{theorem}

\begin{proof} This follows directly from Lemma \ref{Lcf}. The condition $$\sum_{p\leqslant x}f(p)=c\cdot\frac{x}{\log x}(1+o(1))$$ gives by partial summation that \begin{align}\nonumber \sum_{p\leqslant x}\frac{f(p)}{p} &=\frac{1}{x}\sum_{p\leqslant x}f(p)+\int_1^x\left(\sum_{p\leqslant t}f(p)\right)\frac{dt}{t^2}\\
\nonumber &=c\cdot\frac{1}{\log x}(1+o(1))+c\int_1^x\frac{1}{t\log t}(1+o(1))dt\\
\label{cdens}&=c\log\log x(1+o(1)).
\end{align} Note that the proof of the lemma follows from the divergent behavior of $\sum_{p\leqslant x}\frac{f(p)}{p}$ in both \eqref{fposc} and \eqref{fnegc}, and that this divergence is satisfied by \eqref{cdens}. Thus using \eqref{cdens} in the place of \eqref{fposc} and \eqref{fnegc} is enough to prove Lemma \ref{Lcf}, and thus the condition \eqref{cdens} implies the result of the theorem.
\end{proof}

\section{Extension to multiplicative functions}

The results of the previous section are extendable to multiplicative functions $f:\B{N}\to\{-1,1\}$ with the added condition that there is some $k\geq 1$ with $f(2^k)=1$. In this section, by relating a multiplicative function $f:\B{N}\to\{-1,1\}$ to a related completely multiplicative function, we are able to deduce Theorem \ref{Tmain} as a corollary to Theorem \ref{Tcmain}. This is obtained via the following lemma.

\begin{lemma}\label{multf} Let $f:\B{N}\to\{-1,1\}$ be a multiplicative function such that there is some $k\geq 1$ with $f(2^k)=1$. Then $$F(\gs)=\sum_{n\geqslant 1}\frac{f(n)}{n^\gs}=\Pi(\gs)\cdot \prod_p\left(1-\frac{f(p)}{p^\gs}\right)^{-1}\qquad \left(\gs>1\right),$$ where $$\Pi(\gs)=\prod_p\left(1+\sum_{k\geqslant 2}\frac{f(p^k)-f(p^{k-1})f(p)}{p^{k\gs}}\right).$$ Moreover, there is a $\gs_0(f)\in(0,1)$ such that $\Pi(\gs)$ is absolutely convergent for $\gs>\gs_0(f).$
\end{lemma}

\begin{proof} Note that if $f$ is multiplicative, then for $\gs>1$ we have using the Euler product for its generating Dirichlet series that \begin{align*} F(\gs):=\sum_{n\geqslant 1}\frac{f(n)}{n^\gs}&=\prod_p\left(1+\frac{f(p^2)}{p^{2\gs}}+\frac{f(p^3)}{p^{3\gs}}+\cdots\right)\\
&=\prod_p\left(1-\frac{f(p)}{p^{\gs}}\right)^{-1}\cdot\prod_p\left(1+\sum_{k\geqslant 2}\frac{f(p^k)-f(p^{k-1})f(p)}{p^{k\gs}}\right).\end{align*} 

It remains to show that $\Pi(\gs)$ is absolutely convergent for $\gs>\frac{\log \gf}{\log 2}$. Firstly, note that for each prime $p$ we have \begin{multline}\label{Piterm}1+\sum_{k\geqslant 2}\frac{f(p^k)-f(p^{k-1})f(p)}{p^{k\gs}}\\ \geq\max\left\{1+\sum_{k\geqslant 2}\frac{f(2^k)-f(2^{k-1})f(2)}{2^{k\gs}},1-\frac{2}{3^\gs(3^\gs-1)}\right\}.\end{multline} 
To ensure that none of the factors of the product $\Pi(\gs)$ is zero, we will show that the right--hand side of \eqref{Piterm} is greater than zero for $\gs>\frac{\log\gf}{\log 2}.$ Now if $f(2^k)=1$ for all $k\geq 1$, then $$1+\sum_{k\geqslant 2}\frac{f(2^k)-f(2^{k-1})f(2)}{2^{k\gs}}=1>0,$$ regardless of the range of $\gs$. Thus we may suppose that $f(2^k)\neq 1$ identically. Then, using our assumption, since $f(2^k)=1$ for at least one $k\geq 1$, we have for some $k\geq 2$ that $f(2)\neq f(2^k)$. Rephrased, this means that there is a $k\geq 3$ such that $f(2^{k-1})f(2)=-1$. Denote $$k_0:=\min\{k\geq 3: f(2^{k-1})f(2)=-1\}.$$ Then $$1+\sum_{k\geqslant 2}\frac{f(2^k)-f(2^{k-1})f(2)}{2^{k\gs}}\geq 1-2\sum_{k\geq 2}\frac{1}{2^{k\gs}}+\frac{2}{2^{k_0\gs}}=1-\frac{2}{2^\gs(2^\gs-1)}+\frac{2}{2^{k_0\gs}}.$$ Note that for $\gs>0$ and $k_0\geq 2$, the function $$1-\frac{2}{2^\gs(2^\gs-1)}+\frac{2}{2^{k_0\gs}}$$ is continuous and increasing. Also at $\gs=1$ we have $$1-\frac{2}{2^1(2^1-1)}+\frac{2}{2^{k_0}}=\frac{2}{2^{k_0}}>0,$$ so that by continuity and the fact that $k_0\geq 3$, there is some minimal $\ga:=\ga(k_0)\in(0,1)$ such that for $\gs>\ga$ we have $$1-\frac{2}{2^\gs(2^\gs-1)}+\frac{2}{2^{k_0\gs}}>0.$$

Also, we have that $$3^{2\gs}-3^\gs-2>0$$ for $\gs>\frac{\log 2}{\log3}$ by the quadratic formula. Since $3^{2\gs}-3^\gs-2>0$ precisely when $1-\frac{2}{3^\gs(3^\gs-1)}>0$, combining this with the above, we have that each of the terms of the product $\Pi(\gs)$ is positive for all $$\gs>\gs_0(f):=\max\left\{\ga,\frac{\log 2}{\log 3}\right\}.$$ Since this maximum is strictly less that one, the only thing left to show is that the sum $\sum_p\sum_{k\geqslant 2}\frac{f(p^k)-f(p^{k-1})f(p)}{p^{k\gs}}$ is absolutely convergent. We have that \begin{align*} \sum_p\left|\sum_{k\geqslant 2}\frac{f(p^k)-f(p^{k-1})f(p)}{p^{k\gs}}\right| 
\leqslant \sum_p\sum_{k\geqslant 2}\frac{2}{p^{k\sigma}}
=2\sum_p\frac{1}{p^\sigma}\cdot\frac{1}{p^\sigma-1},
\end{align*} which is convergent when $\sigma>1/2$, proving the lemma.
\end{proof}

It is worth remarking that assuming that $f(2^k)=1$ for some $k\geq 1$ ensures that we are not considering the counterexample $g$ defined in \eqref{g}.

We now give the proof of Theorem \ref{Tmain} as a corollary to Theorem \ref{Tcmain}.

\begin{proof}[Proof of Theorem \ref{Tmain}] Let $f:\B{N}\to\{-1,1\}$ be a multiplicative function such that $f(2^k)=1$ for some $k\geq 1$, and denote $F(\gs)=\sum_{n\geqslant 1}\frac{f(n)}{n^\gs}$. By Lemma \ref{multf}, we have that $$F(\gs)=\Pi(\gs)F_c(\gs),$$ where $\Pi(\gs)$ is defined by as in Lemma \ref{multf} and $F_c(\gs)$ is the generating Dirichlet series for the completely multiplicative function $f_c:\B{N}\to\{-1,1\}$ defined by $f_c(p)=f(p)$ for all primes $p$. Similar to \eqref{logF}, we have that $$\log F(\gs)=\log\Pi(\gs)+\log F_c(\gs)=\sum_p\frac{f(p)}{p^\gs}+O(1),$$ since $\Pi(\gs)>0$ for $\gs\geq 1$. If $c>0,$ then considering the proof of Theorem \ref{Tcmain} for $F_c(\gs)$ gives that $$\lim_{\gs\to 1^+}F_c(\gs)=\infty,$$ and so $$\lim_{\gs\to 1^+}F(\gs)=\infty,$$ which in turn gives that the partial sums $\sum_{n\leqslant x} f(n)$ are unbounded. 

If $c<0$ we just consider the proof of Theorem \ref{Tcmain} for the function $1/F(\gs)$, and use the equation $$\log \frac{1}{F(\gs)}=-\log F(\gs)=-\log \Pi(\gs)-\log F_c(\gs)$$ to yield the result.
\end{proof}

\section{Weakening of hypotheses and further extensions}\label{exten}

In Theorem \ref{Tcmain} we can replace the condition \begin{equation}\label{cf}\sum_p\frac{c-f(p)}{p}<\infty\end{equation} with something considerably weaker. 

Note that assumption \eqref{cf} is given so that we may use an asymptotic of the form $$\sum_{p\leqslant x}\frac{f(p)}{p}=c\log\log x+O(1),$$ for nonzero $c\in[-1,1]$. In the case of positive $c$ we can weaken the condition to \begin{equation}\label{liminf}\lim_{x\to\infty}\sum_{p\leqslant x}\frac{f(p)}{p}=\infty,\end{equation} and in the case of negative $c$ we can weaken the condition to \begin{equation}\label{limsup}\lim_{x\to\infty} \sum_{p\leqslant x}\frac{f(p)}{p}=-\infty.\end{equation} Then if \eqref{liminf} holds we have that $\sum_{n\leqslant x} f(n)$ is unbounded, and if \eqref{limsup} holds we have that $\sum_{n\leqslant x} \mu(n)f(n)$ is unbounded. As far as ``density conditions'' the above limits are satisfied when we take $$\sum_{p\leqslant x}f(p) =\frac{c\cdot x}{\log x \log_2 x\cdots \log_k x}(1+o(1)),$$ where $\log_j x$ denotes $\log\log\cdots\log x$ with ``$\log$'' written $j$ times, $k$ is any nonnegative integer, and $c\neq 0$ is taken to be positive or negative depending on the desired case; this is easily seen via partial summation.

We can do a little in the case that $c=0$. Indeed, all we really need is to have for some $\gs>1/2$ that either \begin{equation}\label{weakgs}\lim_{x\to\infty} \sum_{p\leqslant x}\frac{f(p)}{p^\gs}=\infty\qquad\mbox{or}\qquad\lim_{x\to\infty} \sum_{p\leqslant x}\frac{f(p)}{p^\gs}=-\infty.\end{equation} In this case, the same method gives the following theorem, though consideration of the limits in \eqref{weakgs} directly gives a more exact result.

\begin{theorem}\label{extmain} Let $f:\B{N}\to\{-1,1\}$ be a completely multiplicative function, $\gs>1/2$, $k$ a nonnegative integer, and suppose that $$\sum_{p\leqslant x}f(p)=\frac{c\cdot x^\gs}{\log x \log_2 x\cdots \log_k x}(1+o(1)).$$ If $c>0$ then the partial sums of $f$ are unbounded, and if $c<0$ then the partial sums of $\mu f$ are unbounded.
\end{theorem}

As discussed above, the proof of Theorem \ref{extmain} follows exactly the same as that of Theorem \ref{Tcmain}, and as such we omit it for fear of sounding redundant. Theorem \ref{Tmain} can be generalized similarly, but with the added assumptions that both $f\neq g$ for $g$ as defined in \eqref{g}, and that $\gs>\gs_0(f)$ as defined in the proof of Lemma \ref{multf}.

\section{Concluding remarks}

Functions satisfying \eqref{cf} for positive $c$ are, in some sense, large. In fact, since $F(\gs)$ are divergent at $\gs=1$, we have, using an obvious abuse of notation, at least that $$\sum_{n\leqslant x}f(n)\gg x^{1-\eps}$$ for any $\eps>0$. Probably this can be improved, but our original purpose was to just prove the unboundedness of partial sums. Indeed, using the terminology of Goldmakher \cite{Leo}, we should have that the function $f(n)$ mimics the function $c^{\Omega(n)}$ and the partial sums of this function are quite large; we have $$\sum_{n\leqslant x} c^{\Omega(n)}\geqslant \sum_{p\leqslant x}c=c\cdot\pi(x).$$ 

As an extension of the results for negative $c$, it would be nice if one could show that since $\sum_{n\leqslant x}\mu(n)f(n)$ is unbounded, so is $\sum_{n\leqslant x}f(n)$. We suspect that one may have to consider cases whether or not the Riemann hypothesis holds. Nonetheless, since we have $$\sum_{n\leqslant x}\mu(n)f(n)\gg x^{1-\eps}$$ for any $\eps>0$ and the partial sums of $\mu$ are not too small, $\sum_{n\leqslant x} \mu(n)\neq O(x^{1/2}),$ something may be able to be done in this case. Indeed, we conjecture that in this case one should have at least that $$\sum_{n\leqslant x}f(n)\gg x^{1/2-\eps}$$ for any $\eps>0$. Towards something like this we have tried to factor $F(s)$ in an enlightening way (to find a singularity at $s=1/2$, but to no avail. We note that one has for any such series $F(s)$ and $c\in[-1,0)$, that $$F(s)=\left(\frac{\gz(2s)}{\gz(s)}\right)^{|c|}\frac{e^{\frac{P(2s)}{2}}D(s)}{\gz(2s)^{\frac{|c|}{2}}}\cdot\exp\left[-\sum_p\frac{c-f(p)}{p^s}\right],$$ where $P(s)$ is the prime zeta function and the function $D(s)$ is absolutely convergent for $\Re(s)>1/3$.

\bibliographystyle{amsplain}
\providecommand{\bysame}{\leavevmode\hbox to3em{\hrulefill}\thinspace}
\providecommand{\MR}{\relax\ifhmode\unskip\space\fi MR }
\providecommand{\MRhref}[2]{%
  \href{http://www.ams.org/mathscinet-getitem?mr=#1}{#2}
}
\providecommand{\href}[2]{#2}

\end{document}